\documentclass[a4,11pt]{paper}
\usepackage{enumitem,color,amssymb,bm,amsmath,stmaryrd,tikz-cd,mathtools}
\usepackage[english]{babel}  
\usepackage{graphicx}
\usepackage{lipsum}
\usepackage{xcolor}
\usepackage{booktabs}
\usepackage{scalerel}
\usepackage{ntheorem}
{\theoremstyle{empty}
}
\usepackage[top=1in, left=1in, right=1in, bottom=0.8in]{geometry}

\newcommand{\Z}{\mathbb{Z}}

\newcommand{\Line}{{\mathcal{L}}}

\newcommand{\Xvar}{\mathsf{x}}
\newcommand{\Yvar}{\mathsf{y}}
\newcommand{\Tvar}{\mathsf{T}}
\newcommand{\Mvar}{\mathsf{M}}
\newcommand{\Nvar}{\mathsf{N}}

\newcommand{\tpi}{\widetilde{\pi}}
\newcommand{\tB}{\widetilde{B}}
\newcommand{\tz}{\widetilde{z}}
\newcommand{\tX}{\widetilde{X}}

\newcommand{\tw}{\widetilde{w}}

\newcommand{\tDelta}{\widetilde{\Delta}}
\newcommand{\talpha}{\widetilde{\alpha}}
\newcommand{\tbeta}{\widetilde{\beta}}
\newcommand{\tgamma}{\widetilde{\gamma}}
\newcommand{\tGamma}{\widetilde{\Gamma}}
\newcommand{\trho}{\widetilde{\rho}}
\subsectionfont{\large\sf\bfseries\color{black!50!blue}} 
\sectionfont{\Large\sf\bfseries\color{black!50!blue}} 
\titlefont{\LARGE\sf\bfseries\color{black!50!blue}} 
\numberwithin{equation}{section}
\newcommand\p{\mathbf p}

\newcommand{\ra}{\rightarrow}

\newcommand{\Ibb}{\mathbb{I}}

\newcommand{\Qbb}{\mathbb{Q}}
\newcommand{\Rbb}{\mathbb{R}}
\newcommand{\Sbb}{\mathbb{S}}
\newcommand{\Tbb}{\mathbb{T}}

\newcommand{\Zbb}{\mathbb{Z}}
\newcommand{\Cbb}{\mathbb{C}}
\newtheorem{thm}{\bfseries\color{black!50!blue} Theorem}[section]
\newtheorem{prop}[thm]{\bfseries\color{black!50!blue} Proposition}

\newtheorem{lem}[thm]{\bfseries\color{black!50!blue} Lemma}

\newtheorem{rmk}[thm]{\bfseries\color{black!50!blue} Remark}

\newtheorem{defn}[thm]{\bfseries\color{black!50!blue} Definition}

\def\endproof{\relax\ifmmode\expandafter\endproofmath\else
  \unskip\nobreak\hfil\penalty50\hskip.75em\hbox{}\nobreak\hfil\bull
  {\parfillskip=0pt \finalhyphendemerits=0 \bigbreak}\fi}
\def\endproofmath$${\eqno\bull$$\bigbreak}
\def\bull{\vbox{\hrule\hbox{\vrule\kern3pt\vbox{\kern6pt}\kern3pt
\vrule}\hrule}}

\newcommand\x{\mathbf x}
\newcommand\y{\mathbf y}

\newcommand\z{\mathbf z}

\tikzset{commutative diagrams/.cd,
	mysymbol/.style={start anchor=center,end anchor=center,draw=none}
}

\begin{document}
	\title{Holomorphic Triangles in the Heegaard Floer Pillowcase}%
	\author{Eaman Eftekhary and Marzieh Souri}%
	\institution{\sf{School of Mathematics, Institute for Research in Fundamental Sciences (IPM), Tehran, Iran}}
	%————————————————————————————————
	\maketitle
\begin{abstract}
The Heegaard Floer pillowcase is a recurring local model for rational tangle replacement and exact-triangle constructions. We compute the triangle product for its standard rational-slope triples of curves and encode all contributing classes in a finite generating polynomial. One slope ordering gives a vanishing map. For the other, a closed form for the generating function is computed. One may assume that the slopes are $\tfrac10,\tfrac pq$ and $\tfrac{p'}{q'}$.
When an integer lies between the two finite slopes, the closed form has a particularly simple description  and there are exactly $2(pq'-p'q)+2qq'$ contributing classes. We also compute the absolute Maslov bi-gradings of the intersection generators in the employed model, relating the planar labels to the graded Floer groups of two-bridge knots and links.
\end{abstract}

\section{Introduction}

The \emph{pillowcase}, a sphere with four marked points, is a recurring local model in Heegaard Floer theory. Together with related genus-one configurations which may be regarded as its branched covers, it appears in the holomorphic-curve calculations underlying cobordism maps and their compositions. These local calculations connect a broad range of constructions and applications: exact triangles and their unoriented variants \cite{OS-3m2,Ef-triangle,Ciprian-unoriented,AE-1}, bounds on knot invariants such as the $4$-ball genus, unknotting number, and bridge number \cite{OS-four-ball,AE-2,Hom-Wu:tau,AE-unknotting,JMZ}, and obstructions involving tangle replacements, cosmetic operations, and ribbon concordance \cite{Ef-RTR,Ef-torsion,KLMWZ-cosmetic,Zemke-Ribbon}. This recurrence motivates a systematic computation of the triangle product in the pillowcase. For standard rational-slope curves, we encode all contributing classes in a finite generating polynomial, admitting  a simple closed formula whenever an integer lies between the two finite slopes.

The branched double cover of the pillowcase $\Sbb$ is a torus. We use its planar cover to write
\[
\tpi:\mathbb R^2\longrightarrow
\Sbb=\bigl(\mathbb R^2/(2\mathbb Z)^2\bigr)/\{\pm1\}.
\]
The integral lattice projects to $B=\{w_1,w_2,z_1,z_2\}$; even and odd second coordinates specify the $w$- and $z$-types. For a primitive pair $(p,q)$ of integers, let $\alpha_{p/q}$ be the image of $px+qy=1/2$ on the pillowcase $\Sbb$. Set $m=pq'-p'q\neq 0$ and fix the ordered triple
\[
\alpha=\alpha_{1/0},\qquad \beta=\alpha_{p/q}\qquad\text{and}\qquad
\gamma=\alpha_{p'/q'},
\]
where $q,q'>0$ and $p,p'$ are odd. Thus each curve separates the markings into pairs containing one point of each type.
Let $\Delta_{i,j}$ be the planar triangle candidate whose lifted sides lie on
\[
x=\frac12,\qquad px+qy=i+\frac p2\qquad\text{and}\qquad
\qquad p'x+q'y=j+\frac{p'}2.
\]
Triple intersections label infinitesimal candidates, resolved by a small Hamiltonian perturbation. Write $\eta_{i,j}=1$ when the corresponding Maslov-index-zero triangle class $\Delta_{i,j}$ has holomorphic count one modulo two and avoids all markings, and $0$ otherwise. Every contributing class has such a label, with
$\Delta_{i,j}=\Delta_{i+2q,j+2q'}$.
Accordingly, over $\mathbb F=\mathbb F_2$, put
\[
R=\mathbb F[\Xvar,\Yvar]/\langle \Tvar^2-1\rangle,
\qquad\text{where}\qquad \Tvar=\Xvar^q\Yvar^{q'},
\]
and define
\[
\mathcal P_{p/q,p'/q'}(\Xvar,\Yvar)
=\sum_{i=1}^{2q}\sum_{j\in\mathbb Z}\eta_{i,j}\Xvar^i\Yvar^j
=\sum_{j=1}^{2q'}\sum_{i\in\mathbb Z}\eta_{i,j}\Xvar^i\Yvar^j
\in R.
\]
This sum is finite, and both variables are units in $R$.

\begin{thm}\label{thm:1}
	If $p/q<p'/q'$, no Whitney triangle class contributes. In particular, the triangle map for $(\alpha,\beta,\gamma)$ is zero, and therefore $\mathcal P_{p/q,p'/q'}=0$.
\end{thm}

For the opposite ordering, the simplest formula requires only an integer between the slopes.
\begin{thm}\label{thm:2}
	Suppose $m>0$ and
	$\lceil p'/q'\rceil\le\lfloor p/q\rfloor$. Then
	\begin{equation}\label{eq:main-computation}
		\begin{split}
			\mathcal P_{p/q,p'/q'}(\Xvar,\Yvar)
			=\frac{1+\Tvar}{(\Xvar-1)(\Yvar-1)}\Bigl(
			&\Xvar^{\frac{p+1}2}\Yvar^{\frac{p'+1}2}(1-\Xvar^q)+\Xvar^{\frac{1-p}2}\Yvar^{\frac{1-p'}2}(1-\Yvar^{q'})\Bigr).
		\end{split}
	\end{equation}
	There are exactly $2m+2qq'$ contributing Whitney triangle classes.
\end{thm}
The fraction is interpreted in the localization of $R$ at $\Xvar-1$ and $\Yvar-1$. These are non-zero-divisors, and the displayed expression belongs to $R$.

The polynomial determines the triangle map
\[
F_{\alpha\beta\gamma}:C(\alpha,\beta)\otimes C(\beta,\gamma)
\longrightarrow C(\alpha,\gamma),
\]
where $C(\delta,\delta')$ is the Floer complex of the indicated ordered curves, equipped with the differential counting bigons which avoid $B$. These differentials vanish when $\alpha,\beta$ and $\gamma$ are standard, so the intersection generators give homology bases.

The last section treats all slope pairs. In particular, Theorem~\ref{thm:general-row-interval} shows that, for each \(i\), there is an  interval \(I_i=[a_i, b_i]\) with integral endpoints, possibly empty, such that the Whitney triangle $\Delta_{i,j}$ contributes if and only if $j\in I_i$. The endpoints $a_i$ and $b_i$ are given by finite floor and ceiling formulas. Equation~\eqref{eq:general-generating-floor-formula} then determines \(\mathcal P\) without any hypothesis beyond \(m>0\). We prove (see Theorem~\ref{thm:integer-between-slopes}) that the integer-between-slopes condition  \(\lceil p'/q'\rceil\le\lfloor p/q\rfloor\)   is precisely the condition under which the simple generating function  formula (\ref{eq:main-computation}) holds.

The computation fits into the TQFT framework of Ozsv\'ath--Szab\'o \cite{OS-3m1,OS-3m2,OS-knot}, in which holomorphic triangle maps describe cobordisms and polygon relations control their compositions. Explicit calculations on surfaces of low genus also enter the construction and invariance of a link cobordism TQFT \cite{Ef-TQFT}.
From the symplectic perspective, the same triangles compute the binary product $\mu^2$ for the chosen curves in the Fukaya category of the four-punctured sphere, whose $A_\infty$-operations count holomorphic polygons \cite{AS-Fukaya-cat}. Thus the generating polynomial describes a composition law shared by the Heegaard Floer and surface Fukaya settings.

\section{From Heegaard pairs on the sphere to line intersections in the plane}\label{sec:heegaard-triple-const}

\subsection{The double branched cover of the four-punctured sphere}\label{subsec:double-cover}
Throughout this paper, we identify the complex plane $\Cbb$ with $\Rbb^2$. We denote the full integral lattice $\Zbb^2\subset\Rbb^2$ by $\tGamma$ and the sublattice $(2\Zbb)^2$ of even integral points by $\Gamma$. We then set 
\begin{align*}
	T^2:=\Rbb^2/\Zbb^2=\Cbb/\tGamma\qquad\text{and}\qquad
	\Tbb:=\Rbb^2/(2\Zbb)^2=\Cbb/\Gamma.
\end{align*}
Let $\trho:\Cbb\ra T^2$ and $\rho=\rho_{\Gamma}:\Cbb\ra \Tbb$ denote the corresponding quotient map projections. Let $\imath:\Cbb\ra \Cbb$ be the central reflection $\imath(z)=-z$. Since the involution $\imath$ preserves both $\Gamma$ and $\tGamma$ via
\[
\imath (z+\Gamma)=-z+\Gamma\quad\text{and}\quad
\imath (z+\tGamma)=-z+\tGamma
\qquad\forall\ z\in\Cbb,
\]
it induces  well-defined involutions 
\[
\jmath: \Tbb=\Cbb/\Gamma\ra \Cbb/\Gamma=\Tbb\quad\text{and}\quad
\widetilde{\jmath}: T^2=\Cbb/\tGamma\ra \Cbb/\tGamma=T^2.
\]
The quotient of $\Tbb$ under $\jmath$ is a 2-sphere $\Sbb$, and  projection defines a branched double covering 
\[
\pi:\Tbb\ra \Tbb/\jmath=\Sbb.
\]
The branch values of $\pi$ are denoted by $w_1,w_2,z_1,$ and $z_2$. Their pre-images under $\tpi=\pi\circ \rho:\Cbb\ra \Sbb$ form the lattice $\tGamma$. We assign types to these lattice points as follows:

\begin{defn}\label{defn:branched-points}
	 The \emph{type} of a lattice point $(a,b)\in\mathbb{Z}^2$  is defined to be
	 $w_1,w_2,z_1$ or $z_2$ if it is congruent modulo $2$ to 
	 $(0,0),(1,0),(0,1)$ or $(1,1)$, respectively. 
\end{defn}

Lattice points of a given type project under $\tpi$ to the corresponding branch point on $\Sbb$. We set 
\begin{align*}
	B=\{w_1,w_2,z_1,z_2\}\subset \Sbb \quad\text{and}\quad \tB=\{\tw_1,\tw_2,\tz_1,\tz_2\}\subset \Tbb, 
\end{align*}
where $\tw_j$ and $\tz_j$ denote the pre-images of $w_j$ and $z_j$ under $\pi$. Restricting to the complement of these points yields the regular covering maps
\begin{align*}
	\tpi:\Cbb\setminus \tGamma\ra X:=\Sbb\setminus B\quad\text{and}\quad \pi:\tX=\Tbb\setminus \tB\ra X=\Sbb\setminus B.
\end{align*}

\subsection{Rational tangles and simple closed curves on the four-punctured sphere}\label{subsec:single-tangle}
The four-punctured sphere $(\Sbb,B)$ provides the underlying surface for rational tangle diagrams. Realizing $\Sbb$ as the boundary $S^2=\partial D^3$ of the unit 3-ball via
\begin{align*}
	\Sbb=S^2=\{(r,s,t)\in\Rbb^3\mid r^2+s^2+t^2=1\}\quad\text{and}\quad D^3=\{(r,s,t)\in\Rbb^3\mid r^2+s^2+t^2\leq 1\},
\end{align*}    
we place the punctures at
\[
w_1=(-1,0,0),\quad w_2=(1,0,0),\quad z_1=(0,-1,0),\quad\text{and}\quad z_2=(0,1,0).
\]
Each rational tangle corresponds to an isotopy class of simple closed curves $\alpha\subset X=\Sbb\setminus B$ that partitions $B$ into two pairs. Specifically, $\alpha$ separates $\Sbb$ into two disks:
\begin{align*}
	\Sbb\setminus \alpha=D^+_\alpha\amalg D^-_\alpha,\qquad\text{with}\ \ w_1\in D^+_\alpha,
\end{align*}
and represents a rational tangle precisely when 
\[|D^+_\alpha\cap B|=|D^-_\alpha\cap B|=2.\] 
Any such curve bounds a properly embedded disk $(D_\alpha,\alpha =\partial D_\alpha)\subset (D^3, X)$, and connecting the two points in $D^\pm_\alpha\cap B$ by boundary-parallel properly embedded arcs $\gamma_\alpha^\pm \subset D^3\setminus D_\alpha$ yields a tangle $\gamma_\alpha=\gamma^+_\alpha\amalg\gamma^-_\alpha$ that is unique up to boundary-relative isotopy. Thus, rational tangles are in one-to-one correspondence with isotopy classes of essential simple closed curves in $X$ that separate $B$ into 2-element subsets.

For a general closed curve $\alpha\subset X$ (represented by an immersion $\alpha:S^1\ra \Sbb$), the number of connected components of the pullback cover $\alpha^*\tX\to S^1$ is governed by the parity of the branch points enclosed by $\alpha$. Choose a smooth map 
\begin{align*}
	f:(D^2,S^1=\partial D^2)\ra (\Sbb,\alpha(S^1))
\end{align*}
extending $\alpha$ and transverse to $B$. We define
\begin{align*}
	|B\cap\mathrm{Int}(\alpha)|:=\# f^{-1}(B)\pmod2,
\end{align*}
counting points in $f^{-1}(B)$ with multiplicity. This mod-2 count is well-defined and independent of the choice of $f$.

\begin{lem}\label{lem:lift-of-curves}
	Let $\alpha\subset X$ be a closed curve. Then the number $	\#\pi_0(\alpha^*\tX)$ of components of $\alpha^*\tX$ is $2$ if 
	$|B\cap\operatorname{Int}(\alpha)|$ is even, and is $1$ if 
	$|B\cap\operatorname{Int}(\alpha)|$ is odd.
\end{lem}
\begin{proof}
	Since $\pi:\tX\ra X$ is a two-sheeted covering, its monodromy representation
	 \(\nu:\pi_1(X)\ra\mathbb{Z}/2\)
	maps any loop surrounding a single branch point to $1$. If $\alpha$ encloses $k$ branch points (counted with multiplicity), $\nu(\alpha)=k\pmod2$. The pullback two-sheeted cover is therefore disconnected if and only if $k$ is even.
\end{proof}

Consequently, for any curve $\alpha$ representing a rational tangle, $\pi^{-1}(\alpha)$ consists of two disjoint simple closed curves $\alpha^+$ and $\alpha^-$ in $\tX$, interchanged by the involution $\jmath$. The homology class of either lift in the torus determines a primitive direction, up to sign. After isotoping to straight representatives, its image in \(T^2\) has primitive tangent direction \((q,-p)\in \tGamma\), where \(p,q\) are coprime.  
Now consider the family of parallel lines in $\Rbb^2$:
\begin{align*}
	\Line^k_{p/q}=\{(r,s)\in\Rbb^2\mid pr+qs=k\}\subset \Rbb^2\quad \text{for}\ k\in\Rbb.
\end{align*}    
We denote $\Line^{1/2}_{p/q}$ and $\Line^{-1/2}_{p/q}$ by $\Line^+_{p/q}$ and $\Line^-_{p/q}$, respectively. For any lattice point $(r,s)\in\tGamma$, the value $pr+qs$ is an integer; hence, $\Line^k_{p/q}$ intersects $\tGamma$ if and only if $k\in\Zbb$. 

The simple closed image $\trho(\Line^+_{p/q})$ is isotopic to the image of $\alpha^+$ in $T^2$, and $\rho(\Line^+_{p/q})$ is a corresponding representative of $\alpha^+$ on $\Tbb$. Without loss of generality, we take $\alpha^+$ to be $\rho(\Line^+_{p/q})$, which forces $\alpha^-$ to be $\rho(\Line^-_{p/q})$. We call this the \emph{standard representation} $\alpha_{p/q}$. Note that $\alpha_{(-p)/(-q)}=\alpha_{p/q}$; thus, the isotopy class of $\alpha$ is uniquely specified by the slope $\tfrac{p}{q}\in\Qbb\cup\{\infty\}$, where $\infty$ corresponds to $(p,q)=(1,0)$. For non-infinite slopes, we adopt the convention $q>0$.

The pre-images of $\alpha_{p/q}$ under $\tpi:\Rbb^2\ra \Sbb$ are the lines $\Line_{p/q}^{k+1/2}$ for $k\in\Zbb$. These lines partition $\Rbb^2$ into parallel strips $S^k_{p/q}$ bounded by $\Line_{p/q}^{k-1/2}$ and $\Line_{p/q}^{k+1/2}$. We color a strip white if it contains lattice points of type $w_1$, and gray otherwise. A rational tangle is called \emph{orientation-reversing} if its white strips contain lattice points of types $w_1$ and $w_2$, and \emph{orientation-preserving} otherwise. These are related to different notions of rational replacements, as discussed in \cite{Ef-RTR}.

\begin{lem}
	A rational tangle $\alpha_{p/q}$ is orientation-preserving if and only if $p$ is odd.
\end{lem}    
\begin{proof}
	Every lattice point $(r,s)\in\tGamma$ lies in the strip $S^{pr+qs}_{p/q}$. If $p$ is even, $q$ must be odd, so $pr+qs\equiv s\pmod 2$. In this case, points of types $w_1$ and $w_2$ lie in $S^k_{p/q}$ for even $k$, whereas points of types $z_1$ and $z_2$ lie in $S^k_{p/q}$ for odd $k$, making $\alpha_{p/q}$ orientation-reversing.
	
	If $p$ is odd and $q$ is even, $pr+qs\equiv r\pmod 2$. Here, points of types $w_1$ and $z_1$ lie in white strips while $w_2$ and $z_2$ lie in gray strips, so $\alpha_{p/q}$ is orientation-preserving.
	
	Finally, if both $p$ and $q$ are odd, $pr+qs\equiv r+s\pmod 2$. Thus $w_1$ and $z_2$ lie in white strips while $w_2$ and $z_1$ lie in gray strips, again making $\alpha_{p/q}$ orientation-preserving.
\end{proof}    

\subsection{Pairs of rational tangles and line intersections}\label{subse:pairs-of-tangles}
Let $\alpha=\alpha_{p/q}$ and $\beta=\alpha_{p'/q'}$ correspond to distinct rational tangles with distinct slopes in $\Qbb\cup\{\infty\}$. We use primitive pairs with $q,q'\ge0$, representing $\infty$ by $(1,0)$ and define
\begin{align*}
	\sigma=\sigma_{\frac{p}{q},\frac{p'}{q'}}=
	\begin{pmatrix}p&q\\ p'&q'\end{pmatrix}\quad\text{and}\quad
	m=m_{\frac{p}{q},\frac{p'}{q'}}=\det(\sigma)=pq'-p'q.
\end{align*}    
The pre-images of $\alpha$ and $\beta$ under $\tpi:\Rbb^2\ra \Sbb$ are two families of parallel lines:
\begin{align*}
	\tpi^{-1}(\alpha)=\bigcup_{k\in\Zbb}\Line^{k+\frac{1}{2}}_{p/q}\qquad\text{and}\qquad
	\tpi^{-1}(\beta)=\bigcup_{k\in\Zbb}\Line^{k+\frac{1}{2}}_{p'/q'}.
\end{align*}    
Their intersections $p_{k,l}=p_{k,l}(\tfrac{p}{q},\tfrac{p'}{q'})\in\Rbb^2$ are determined by the linear systems of two equations
\begin{align*}
&pr+qs=k+\frac{1}{2}&&\text{and} &&p'r+q's=l+\frac{1}{2}.
\end{align*}

\begin{lem}\label{lem:line-intersection-coords}
	For distinct slopes $\frac{p}{q},\frac{p'}{q'}\in\Qbb\cup\{\infty\}$, setting $m=pq'-p'q$, the lines $\Line^{k+1/2}_{p/q}$ and $\Line^{l+1/2}_{p'/q'}$ intersect at
	\begin{align*}
		p_{k,l}=p_{k,l}(\tfrac{p}{q},\tfrac{p'}{q'})
		&=\Big(\frac{q'(2k+1)-q(2l+1)}{2m},\frac{p(2l+1)-p'(2k+1)}{2m}\Big).
	\end{align*}    
\end{lem}    

Two intersection points $p_{k,l}$ and $p_{k',l'}$ project to the same point on $\Tbb$ if and only if 
\[p_{k,l} - p_{k',l'} \in \Gamma = (2\Zbb)^2,\] 
which holds when
\begin{align*}
	2im=q'(k-k')-q(l-l')\qquad\text{and}\qquad 2jm=p(l-l')-p'(k-k')
\end{align*}    
for some integers $i,j\in\Zbb$. Equivalently, this translates to
\begin{align*}
	(k,l)-(k',l')=2i(p,p')+2j(q,q').
\end{align*}    
Thus, the intersection points of $\pi^{-1}(\alpha)$ and $\pi^{-1}(\beta)$ on $\Tbb$ are parameterized by 
\begin{align*}
	\Ibb_{\frac{p}{q},\frac{p'}{q'}}=\frac{\Zbb\oplus\Zbb}{\langle(2p,2p'),(2q,2q')\rangle}=\frac{\Zbb\oplus\Zbb}{\Gamma_{\frac{p}{q},\frac{p'}{q'}}},
\end{align*}    
where $\Gamma_{p/q,p'/q'} = \sigma(\Gamma)$. Since $\det(\sigma)=m$, the finite group $\Ibb_{p/q,p'/q'}$ has order $4|m|$, corresponding to $4|m|$ intersection points on $\Tbb$.

The central reflection $(r,s)\mapsto (-r,-s)$ on $\Rbb^2$ induces the involution $\tau(k,l)=(-k-1,-l-1)$ on the grid parameters $(k,l)$. Note that
\begin{align*}
	\tau\Big((k,l)+\Gamma_{\frac{p}{q},\frac{p'}{q'}}\Big)
	=\tau(k,l)+\Gamma_{\frac{p}{q},\frac{p'}{q'}},
\end{align*}    
so $\tau$ induces a well-defined action on $\Ibb_{p/q,p'/q'}$. The quotient $\Ibb_{p/q,p'/q'}/\langle \tau \rangle$ consists of $2|m|$ orbits, which are in unique correspondence with the intersection points of $\alpha$ and $\beta$ downstairs on $X$. We summarize this in the following lemma:

\begin{lem}\label{lem:line-intersection}
	Given distinct rational slopes $\tfrac{p}{q},\tfrac{p'}{q'}\in\Qbb\cup\{\infty\}$, the intersection points between the curves $\alpha_{p/q}$ and $\alpha_{p'/q'}$ on $X$ are in natural bijection with 
	\begin{align*}
		\bigl(\Zbb^2/\sigma(\Gamma)\bigr)/\langle\tau\rangle,
	\end{align*}    
	where $\sigma=\begin{pmatrix}p&q\\p'&q'\end{pmatrix}$ and $\tau:\Zbb^2\ra \Zbb^2$ is the involution $\tau(k,l)=(-k-1,-l-1)$.
\end{lem}

\subsection{Maslov bi-gradings}\label{subsec:line-maslov-bigrading}
We retain the notation
\(p_{k,l}=\Line_{p/q}^{k+1/2}\cap \Line_{p'/q'}^{l+1/2}\) of Lemma~\ref{lem:line-intersection-coords}, and write
\(\operatorname{gr}_{\alpha,\beta}(p_{k,l})\) for the Maslov bi-grading of its image under \(\tpi\), with the attaching curves ordered as \((\alpha,\beta)\).
Assume that both pairs are primitive, \(p,p'\) are odd, and \(m=pq'-p'q>0\).
The ordered pair specifies a four-pointed genus-zero Heegaard diagram for a two-bridge knot or link in \(S^3\), oriented by the \(w/z\)-markings.
Positivity can be achieved by interchanging the curves, although reversing their order replaces \((M_w,M_z)\) by \((-1-M_w,-1-M_z)\).
With the surface orientation induced from the plane, our convention is
\begin{align*}
	M_w(\x)-M_w(\y)=\mu(D)-2n_w(D),\qquad
	M_z(\x)-M_z(\y)=\mu(D)-2n_z(D),
\end{align*}
where the positive \(\alpha\)-boundary of \(D\) runs from \(\x\) to \(\y\), and multiplicities are summed over each type.
We normalize each Maslov grading so that the homology after forgetting the markings of the other type is supported in degrees \(0\) and \(-1\).

For \(m>0\), and \(r,t\in\Zbb\) with \(\gcd(2r+1,m)=1\), define
\begin{align*}
	F_{m,r}(t,n)
	=-\sum_{i=t}^{n+t-1}(-1)^{\left\lfloor\frac{(2r+1)i}{m}\right\rfloor}
	\qquad(n\geq 0),
\end{align*}
and set \(F_{m,r}(t,-n)=-F_{m,r}(t-n,n)\) when \(n>0\).
Empty sums are zero. The summands change sign under \(i\mapsto i+m\), so the sum over any \(2m\) consecutive integers vanishes.
Consequently, \(F_{m,r}(t,n)\) depends only on the residues of $t$ and $n$  modulo \(2m\).

\begin{thm}\label{thm:line-maslov-bigrading}
	Let \(\tfrac pq,\tfrac{p'}{q'}\), \(m=pq'-p'q>0\), and \(p_{k,l}\) be as above.
	Choose \(a,b\in\Zbb\) with \(pa+qb=1\) and \(b\) even, and set
	\(2r+1=p'a+q'b\).
	Choose \(s\in\{0,\ldots,m-1\}\) so that \(2s+1\) is the inverse of \(2r+1\) modulo \(2m\), and put
	\begin{align}\label{eq:value-of-n}
		n=(-1)^k\bigl(l-(2r+1)k-r\bigr).
	\end{align}
	Then
	\begin{align*}
		\operatorname{gr}_{\alpha,\beta}\bigl(p_{k,l}(\tfrac pq,\tfrac{p'}{q'})\bigr)
		&=\frac{F_{m,r}(0,m)}2(1,1)
		+\left(F_{m,s}(-r,n)-\frac{F_{m,s}(-r,m)}2\right)(1,-1).
	\end{align*}
	The result is integral, independent of the B\'ezout choice, and invariant under the deck identifications of Lemma~\ref{lem:line-intersection}.
\end{thm}
\begin{proof}
    A B\'ezout solution with \(b\) even exists because replacing \((a,b)\) by \((a+q,b-p)\) changes the parity of \(b\).
	The coordinate change \(u=px+qy,\ v=-bx+ay\) has determinant one and preserves the deck group and the \(w/z\)-types, since
	\(y=bu+pv\equiv v\pmod2\) on the lattice.
	Since the coordinate change is unimodular, \((2r+1,m)\) is primitive. Moreover, \(a\) is odd, so \(2r+1=p'a+q'b\) is odd. Hence \(\gcd(2r+1,2m)=1\), ensuring that \(2s+1\) exists.
	The selected lines become
	\begin{align*}
		u=k+\tfrac12,\qquad (2r+1)u+mv=l+\tfrac12.
	\end{align*}
	Translate by \((-k,0)\) if \(k\) is even, or apply the half-turn
	\((u,v)\mapsto(-u+k+1,-v)\) if \(k\) is odd.
	These deck transformations send the intersection to the point \((\tfrac12,\tfrac nm)\) on line $x=\tfrac12$, where $n$ is determined by (\ref{eq:value-of-n}).
	Thus it suffices to grade the intersection points for the reduced pair \((\alpha_{1/0},\alpha_{(2r+1)/m})\).
	
	Let \(\x_i\) be the image of \((\tfrac12,\tfrac im)\).
	Choose a lattice center \((u_0,v_0)\) satisfying
	\begin{align*}
		(2r+1)u_0+mv_0=i+r+1\qquad\text{and}\qquad \big|u_0-\tfrac12\big|\le\tfrac m2.
	\end{align*}
	The first coordinates of the integral solutions of the first equation form one residue class modulo \(m\), so a nearest representative to \(\tfrac12\) works.
	Let \(Q_i\) be the parallelogram bounded by
	\begin{align*}
		u=\tfrac12,\quad u=2u_0-\tfrac12,\qquad
		(2r+1)u+mv=i+r+\tfrac {1}2,\quad (2r+1)u+mv=i+r+\tfrac {3}2.
	\end{align*}
	Two vertices lift to \(\x_i\) and \(\x_{i+1}\); the other two are their half-turn images about \((u_0,v_0)\).
	There are no boundary lattice points, since the defining levels are half-integral.
	Every interior lattice point lies on \((2r+1)u+mv=i+r+1\), and its first coordinate differs from \(u_0\) by a multiple of \(m\).
	The distance of its first coordinate from \(u_0\) is less than \(|u_0-\tfrac12|\le \tfrac m2\), so the center is the only lattice point in \(Q_i\).
	
	The half-turn quotient of \(Q_i\) is a disk with two convex corners, whose image on the pillowcase is an immersed index-one bigon.
	It contains one marked point with multiplicity one, of type \(w\) if \(v_0\) is even and type \(z\) otherwise.
	For \(u_0\le0\), its positive \(\alpha\)-boundary runs from \(\x_i\) to \(\x_{i+1}\); for \(u_0\ge1\), the direction is reversed.
	The relative grading convention therefore gives
	\begin{align*}
		M_w(\x_{i+1})-M_w(\x_i)
		&=\operatorname{sgn}(1-2u_0)(-1)^{v_0}
		=-(-1)^{\left\lfloor\frac{(2s+1)(i+r+1)-1}{m}\right\rfloor}
		=-(-1)^{\left\lfloor\frac{(2s+1)(i-r)}{m}\right\rfloor},
	\end{align*}
	and the \(M_z\)-increment is its negative.
	For the middle equality, reduce \((2s+1)(i+r+1)\) modulo \(2m\): the sign is negative on \(\{1,\ldots,m\}\) and positive on its complement.
	The center equation gives
	\[(2s+1)(i+r+1)\equiv u_0+mv_0\pmod{2m},\] 
	and the chosen bound on \(u_0\) gives precisely these signs. For the last equality, the two numerators differ by
	\((2s+1)(2r+1)-1\), a multiple of \(2m\);
	hence their floors differ by an even integer.
	
	Summing the increments gives, for every integer \(n\),
	\begin{align*}
		M_w(\x_n)=c_w+F_{m,s}(-r,n),\qquad
		M_z(\x_n)=c_z-F_{m,s}(-r,n),
	\end{align*}
	for some constants $c_w$ and $c_z$.
	Antiperiodicity of the summands implies
	\begin{align*}
		F_{m,s}(-r,n+m)
		=F_{m,s}(-r,m)-F_{m,s}(-r,n).
	\end{align*}
	Translation by \((0,1)\) sends \(\x_n\) to \(\x_{n+m}\)
	and interchanges the basepoint types. It therefore exchanges
	the two normalized gradings, giving \(c_z-c_w=F_{m,s}(-r,m)\).

	Let \(L\) be the oriented two-bridge link of the reduced ordered pair. The Schubert signature formula \cite[Theorem~9.3.6]{Murasugi} is
	\begin{align*}
		\sigma(L)=-\sum_{i=1}^{m-1}(-1)^{\lfloor \frac{(2r+1)i}{m}\rfloor}
		=F_{m,r}(0,m)+1,
	\end{align*}
	with the convention that \((\tfrac10,\tfrac13)\) has signature \(-2\). The numerator \(2r+1\) may first be reduced modulo \(2m\), which leaves the sum unchanged.
	Alternating-link thinness in our  normalization gives
	\begin{align*}
		A=\frac{M_w-M_z-1}{2},\qquad M_w-A=\frac{\sigma(L)}2;
	\end{align*}
	here \(A\) is the collapsed Alexander grading, with the normalization and extra-basepoint factors as in \cite{Wong-2017}. 
	The differential avoiding all four marked points vanishes: a contributing bigon would lift to a bigon between two straight lines, which is impossible.
	Thus thinness applies to every generator, and
	\(c_w+c_z=\sigma(L)-1=F_{m,r}(0,m)\).
	Solving for the constants proves the formula.
	Both length-\(m\) sums are congruent to \(m\pmod2\), proving integrality.
	Finally, changes of lift or B\'ezout coordinates preserve the marked, ordered diagram and hence its normalized gradings.
\end{proof}

\section{Triples of Rational Tangles and Whitney Triangles}\label{sec:triangle-labels}
Let three orientation-preserving simple closed curves $\alpha,\beta,\gamma\subset \Sbb\setminus B$ specify three non-isotopic rational tangles, corresponding to three distinct extended rational numbers in $\Qbb\cup\{\infty\}$. By applying a marked integral coordinate change as in the proof of Theorem~\ref{thm:line-maslov-bigrading}, we fix $\alpha$ to correspond to $\infty=\tfrac{1}{0}$, while $\beta=\alpha_{p/q}$ and $\gamma=\alpha_{p'/q'}$ correspond to finite rational numbers with $q,q'>0$. We further assume that $p$ and $p'$ are odd, ensuring that all three tangles are orientation-preserving. We also set $m=pq'-p'q$. For the remainder of the paper, these parameters will be fixed.

\subsection{Triangles in the planar cover}
Given intersection points $\x\in\alpha\cap\beta$, $\y\in\beta\cap\gamma$, and $\z\in\gamma\cap\alpha$, the set of homotopy classes of Whitney triangles with vertices at $(\x,\y,\z)$ and boundary arcs along $(\alpha,\beta,\gamma)$ is denoted by $\pi_2(\x,\y,\z)$. In the domain of the standard topological triangle, we denote the vertices mapped to $\x,\y,$ and $\z$ by $v_\gamma$, $v_\alpha$, and $v_\beta$, respectively, and label the corresponding boundary arcs as the $\alpha$-side, $\beta$-side, and $\gamma$-side.
We denote the set of all Whitney triangles in the Heegaard triple diagram $(\Sbb,\alpha,\beta,\gamma,B)$ by
\begin{align*}
	\pi_2(\alpha,\beta,\gamma):=\coprod_{\x\in\alpha\cap\beta}\coprod_{\y\in\beta\cap\gamma}\coprod_{\z\in\gamma\cap\alpha}\pi_2(\x,\y,\z). 
\end{align*}
For a class $\Delta\in\pi_2(\alpha,\beta,\gamma)$, we set $\eta(\Delta)=1$ if $\Delta$ has Maslov index zero, its multiplicity at every point of $B$ is zero, and its moduli space of holomorphic representatives has count $1\in\Zbb/2$; otherwise, we set $\eta(\Delta)=0$. This yields the indicator map
\begin{align*}
	\eta:\pi_2(\alpha,\beta,\gamma)\ra \Zbb/2,
\end{align*}
and a partition
\begin{align*}
	\pi_2(\alpha,\beta,\gamma)=\pi_2^\circ(\alpha,\beta,\gamma)\amalg\pi_2^{\bullet}(\alpha,\beta,\gamma),\qquad
	\text{where}\ \ \pi_2^\bullet(\alpha,\beta,\gamma):=\{\Delta\in\pi_2(\alpha,\beta,\gamma)\mid \eta(\Delta)=1\}.
\end{align*}
We call the classes in $\pi_2^{\bullet}(\alpha,\beta,\gamma)$ the \emph{contributing Whitney triangles}. A Whitney triangle whose image avoids $B$ lifts through $\tpi:\Rbb^2\setminus\tGamma\ra\Sbb\setminus B$.  Conversely, triangular regions in $\Rbb^2$ project to planar candidate classes via the branched cover $\tpi=\pi\circ\rho:\mathbb{R}^{2}\ra \Sbb$.

As in Section~\ref{sec:heegaard-triple-const}, we define the line families in $\Rbb^2$ for all $k\in\Zbb$:
\begin{align*}
	\Line_k^\alpha=\Line^{k+\frac{1}{2}}_{1/0},\quad \Line_k^\beta=\Line^{k+\frac{1}{2}}_{p/q}, \quad\text{and}\quad \Line_k^\gamma=\Line^{k+\frac{1}{2}}_{p'/q'},
\end{align*}
and set
\begin{align*}
	\widetilde{\alpha}=\left\{\Line_k^\alpha\right\}_{k\in\Zbb},\quad
	\widetilde{\beta}=\left\{\Line_k^\beta\right\}_{k\in\Zbb},
	\quad\text{and}\quad\widetilde{\gamma}=\left\{\Line_k^\gamma\right\}_{k\in\Zbb}.
\end{align*}

When the lines $\Line_k^\alpha,\Line_l^\beta$ and $\Line_m^\gamma$ all pass through the same point, we obtain an {\emph{degenerate}} triangle. In this case, we need to perturb the picture by a Hamiltonian isotopy which is applied to $\alpha$ in order to make sure that there is no concurrency of three lines. 
The degenerate triangle then becomes a very small embedded triangle away from all lattice points. For non-degenerate triangles (those that are not degenerate) the small Hamiltonian isotopy may simply be ignored.

The following definition specifies the lifted regions corresponding to contributing Whitney triangles:

\begin{defn}\label{defn:accebtable-triangles}
	A non-degenerate triangular region $\tDelta\subset\mathbb{R}^{2}$ is an \emph{acceptable triangle} if it satisfies:
	\begin{enumerate}
		\item \textbf{Lifting:} Each of its three sides lies on a line component of $\widetilde{\alpha}$, $\widetilde{\beta}$, or $\widetilde{\gamma}$ (termed its $\alpha$-side, $\beta$-side, and $\gamma$-side, respectively).
		\item \textbf{Positivity:} The $\alpha$-side, $\beta$-side, and $\gamma$-side occur in clockwise order along $\partial\tDelta$.
		\item \textbf{Emptiness:} The triangle contains no integral lattice points: $\tDelta\cap\tGamma=\varnothing$. 
	\end{enumerate}
	An acceptable triangle $\tDelta$ is \emph{$\alpha$-normal} if its $\alpha$-side lies on $\Line_0^\alpha$. ($\beta$-normal and $\gamma$-normal triangles are defined analogously.)
\end{defn}

The following proposition treats non-degenerate triangles in the straight-line model.

\begin{prop}\label{prop:whitney-triangles}
	For non-degenerate candidates in the straight-line model, every contributing Whitney triangle $\Delta\in\pi_2^{\bullet}(\alpha,\beta,\gamma)$ is the projection under 
	\[\tpi=\pi\circ\rho:\Rbb^2\setminus\tGamma\ra \Sbb\setminus B\] 
	of an acceptable triangle $\tDelta\subset\Rbb^2$ with boundaries on $\talpha,\tbeta,$ and $\tgamma$.
	Conversely, every acceptable triangle projects to a contributing class. Furthermore, $\tDelta$ can be chosen to be $\alpha$-normal. Any other $\alpha$-normal acceptable triangle $\tDelta'$ projecting to $\Delta$ is of the form \(\tDelta'=\tDelta+(0,2k)\) for some \(k\in\Zbb\).
\end{prop}
\begin{proof}
A triangle avoiding \(B\) lifts through the unbranched covering \(\Rbb^2\setminus\tGamma\to\Sbb\setminus B\). The lift is unique once a corner is chosen, and each boundary arc lies on one fixed line. The local surface polygon correspondence identifies rigid triangles with immersed disks having three convex corners; see \cite[Appendix B, Proposition B.1]{Azam-Blanchet}. 
Each lifted boundary arc is an immersion into a line and hence is monotone. The three arcs therefore form the boundary of a Euclidean triangle. Their winding number is one in its interior and zero outside, so the lifted holomorphic map has degree one onto this triangle. Positivity gives the required order of its sides, and avoidance of \(B\) gives emptiness.

Conversely, the Riemann mapping theorem parametrizes a positive, lattice-free Euclidean triangle uniquely, with its three corners prescribed. Projection gives an immersed holomorphic triangle on the pillowcase. Since the covering is a local biholomorphism along the entire image, it identifies the pulled-back tangent bundles and boundary conditions, preserving the Fredholm index. The surface index formula is $\mu=2e-\tfrac12$ for a triangle; here $e=\tfrac14$, so $\mu=0$ \cite[Appendix B]{Azam-Blanchet}. Its count is therefore one modulo two. This computation uses the immersed source disk and does not require its projected image to be embedded. 

Deck transformations place the vertical side on \(\Line_0^\alpha\). The stabilizer of this line in the deck group consists exactly of translations by \((0,2k)\), proving the final assertion.
\end{proof}

Having settled non-degenerate triangles, we now consider triple intersections.
A sufficiently small Hamiltonian perturbation of $\alpha$, fixing $\beta,\gamma$, preserves the contribution criterion for non-degenerate planar candidates and resolves each triple intersection into one small triangle. It is easy to see that the small triangle has correct orientation (and therefore contributes) exactly when $m=pq'-p'q>0$. Pairwise intersections retain their labels when we do the perturbations.

From now on, triangle counts refer to this perturbed diagram and planar coordinates refer to its straight-line limit. Pairwise differentials still vanish, since each pair of lifted curves intersects at most once and cannot bound a bigon. In particular, the degenerate case may be treated just as the non-degenerate case.

\begin{rmk}
	The line $\Line^\alpha_0$ is given by $\{(\tfrac{1}{2},s)\mid s\in\Rbb\}$. For any planar candidate $\Delta$, we uniquely fix its \emph{standard lift} $\tDelta$ by requiring $\tDelta(v_\gamma)=(\tfrac{1}{2},s)$ with $s\in (0,2]$. Every contributing class possesses a unique standard lift.
\end{rmk}

Let $\tDelta$ be a standard planar candidate. Since $\tDelta(v_\gamma)\in \Line^\beta_k$, Lemma~\ref{lem:line-intersection-coords} gives
\begin{align*}
	\tDelta(v_\gamma)=p_{k,0}(\tfrac{p}{q},\tfrac{1}{0})=\Big(\frac{1}{2},\frac{2k+1-p}{2q}\Big).
\end{align*}
The standard normalization $s=\tfrac{2k+1-p}{2q}\in (0,2]$ with $q>0$ requires
\begin{align*}
	k\in \Big(\frac{p-1}{2}, \frac{p-1}{2}+2q \Big]. 
\end{align*}
The $2q$ integer choices for $k$ in this interval correspond to the $2q$ intersection points in $\alpha\cap\beta$. Since $\beta$ is orientation-preserving, $p$ is odd, making $\tfrac{p-1}{2}$ an integer. We may thus parameterize $k$ as $k=\frac{p-1}{2}+i$ for $i\in\{1,\ldots,2q\}$, yielding
\begin{align*}
	p_{k,0}(\tfrac{p}{q},\tfrac{1}{0})=\Big(\frac{1}{2},\frac{i}{q}\Big)=:p^\gamma_i.
\end{align*}

Similarly, for the vertex $v_\beta$ (the $\alpha\cap\gamma$ intersection), $\tDelta(v_\beta)$ has the form
\begin{align*}
	p_{l,0}(\tfrac{p'}{q'},\tfrac{1}{0})=\Big(\frac{1}{2},\frac{2l+1-p'}{2q'}\Big)
	=\Big(\frac{1}{2},\frac{j}{q'}\Big)=:p^\beta_j, \qquad\text{where}\ \ \ l=\frac{p'-1}{2}+j. 
\end{align*} 

Finally, with $k$ and $l$ determined by $i$ and $j$, the lifted image of the third vertex $v_\alpha$ (the $\beta\cap\gamma$ intersection) is 
\begin{align*}
	p_{k,l}(\tfrac{p}{q},\tfrac{p'}{q'})
	&=\Big(\frac{q'(2k+1)-q(2l+1)}{2m},\frac{p(2l+1)-p'(2k+1)}{2m}\Big)\\
	&=\Big(\frac{1}{2}+\frac{q'i-qj}{m},\frac{pj-p'i}{m}\Big)=:p^\alpha_{i,j}.
\end{align*}

We summarize this parametrization in the following proposition:

\begin{prop}\label{prop:vertices-of-triangle}
	The standard planar candidates are parameterized by pairs of integers $(i,j)$ with $i\in\{1,\ldots,2q\}$ and $j\in\Zbb$ and are denoted by $\tDelta_{i,j}$. Every contributing class occurs among these candidates. The vertices of $\tDelta_{i,j}$ are
	\begin{align*}
		p^\gamma_i=\Big(\frac{1}{2},\frac{i}{q}\Big),\quad 
		p^\beta_j=\Big(\frac{1}{2},\frac{j}{q'}\Big),\quad\text{and}\quad
		p^\alpha_{i,j}=\Big(\frac{1}{2}+\frac{q'i-qj}{m},\frac{pj-p'i}{m}\Big).
	\end{align*}    
\end{prop}    

The bi-gradings of the vertices follow from Theorem~\ref{thm:line-maslov-bigrading}. For \(p_i^\gamma\) use the ordered pair \((\alpha,\beta)\) and line indices \((0,(p-1)/2+i)\); for \(p_j^\beta\) use \((\alpha,\gamma)\) and \((0,(p'-1)/2+j)\). For \(p_{i,j}^\alpha\) use \((\beta,\gamma)\) and
\[
 (k,l)=\Big(\frac{p-1}{2}+i,\frac{p'-1}{2}+j\Big).
\]
Note that the output order is \((\alpha,\gamma)\), rather than its reverse.

\section{Classification of the Lifted Triangles}\label{sec:main-triangles}

\subsection{The main theorem}
We denote the planar candidate corresponding to $(i,j)\in\{1,\ldots,2q\}\times\Zbb$ by $\Delta_{i,j}$ and its standard lift by $\tDelta_{i,j}$, defined using the lines $\Line^\alpha_0$, $\Line^\beta_{k_i}$, and $\Line^\gamma_{l_j}$ with 
\[k_i=\frac{p-1}{2}+i\qquad\text{and}\qquad l_j=\frac{p'-1}{2}+j.\]

Two of the vertices of $\tDelta_{i,j}$ are on the vertical line $x=\tfrac12$. The horizontal displacement of the third vertex is $(q'i-qj)/m$. If we set $D=q'i-qj$, the sign of $D$ also indicates whether $p_i^\gamma$ is above $p_j^\beta$ (the case $D>0$), $p_i^\gamma$ is below $p_j^\beta$ (the case $D<0$), or they coincide (the degenerate case $D=0$).

\begin{prop}\label{prop:orientation-test}
	If $\tfrac{p}{q}>\tfrac{p'}{q'}$, all 
	planar triangles $\Delta_{i,j}$ (for $i\in\{1,\ldots,2q\}$ and $j\in\Zbb$) are positively oriented. If $\tfrac{p}{q}<\tfrac{p'}{q'}$, all  planar triangles $\Delta_{i,j}$ are negatively oriented.
\end{prop}    
\begin{proof}
	As before, set \(D=q'i-qj\). If $D=0$, the corresponding triangle is degenerate and we have already established the claim. For $D\neq 0$,
	the signed area in the vertex order
	\((p_i^\gamma,p_{i,j}^\alpha,p_j^\beta)\) is
	\[
	\frac12\det\bigl(p_{i,j}^\alpha-p_i^\gamma,
	p_j^\beta-p_i^\gamma\bigr)
	=-\frac{D^2}{2mqq'}.
	\]
	The order of vertices  is clockwise exactly when \(m>0\), which is the prescribed positive side order. 
\end{proof}

Let us thus assume that $m>0$.
 For $i\in\Zbb$ and $\ell\in\Zbb_{\ge0}$, define the lower bounds $a_i^{\ell}$ and the upper bounds $b_i^{\ell}$ by 
\begin{align*}
	a_i^{\ell}&=\frac{p'+1}{2}+p'\ell
	+q'\Big\lfloor\frac{2i-p(2\ell+1)}{2q}\Big\rfloor\qquad\text{and}\qquad
	b_i^{\ell}=-\frac{p'+1}{2}-p'\ell
	+q'\Big\lceil\frac{2i+p(2\ell+1)}{2q}\Big\rceil.
\end{align*}
Then $b_i^{\ell}=-a_{-i}^{\ell}$. Define
\begin{equation}\label{eq:general-row-endpoints}
	 a_i=\max_{0\le\ell<q}a_i^{\ell}
	\qquad\text{and}\qquad
	 b_i=\min_{0\le\ell<q}b_i^{\ell}
	=- a_{-i}.
\end{equation}

\begin{thm}\label{thm:general-row-interval}
	If $m<0$, no Whitney triangle class contributes.
	If $m>0$, the candidate $\Delta_{i,j}$ contributes if and only if $j$ is in the interval $I_i=[a_i,b_i]$ with integer endpoints.
	Moreover, we have
	\[ a_{i+q}= a_i+q',\qquad
	 b_{i+q}= b_i+q',\qquad
	 b_i=- a_{-i}\qquad\text{and}\qquad
	 b_i- a_i+1\ge0.
	\]
\end{thm}
\begin{proof}
    The assertion for $m<0$ follows from Proposition~\ref{prop:orientation-test}. Assume $m>0$; Proposition~\ref{prop:whitney-triangles} 
    reduces the question to lattice avoidance, including the degenerate case.
	Fix $i$, and put $c=q'i/q$. 
	We first test all integer columns to the right of $x=\tfrac12$, including columns beyond the third vertex. At $x=\ell+1$, with $\ell\geq 0$ an integer, the two affine side heights are
	\[
	h_\ell=\frac{i-p(\ell+1/2)}q,
	\qquad g_\ell=\frac{j-p'(\ell+1/2)}{q'}.
	\]
	If $j<c$, the third vertex of the triangle is on the right of $x=\tfrac12$. When the  column $x=\ell+1$ is reached, its slice cut by the triangle is $[g_\ell,h_\ell]$, which avoids integers exactly when \(g_\ell>\lfloor h_\ell\rfloor\).
	If the column is beyond the vertex, then $g_\ell>h_\ell$, and the same inequality holds automatically. The condition \(g_\ell>\lfloor h_\ell\rfloor\) is equivalent to
	\[
	j>p'\Big(\ell+\frac12\Big)+q'\lfloor h_\ell\rfloor,
	\]
	and, since $p'$ is odd and $j$ is integral, to $j\ge a_i^{\ell}$.
	
	For $j>c$, the third vertex of $\Delta_{i,j}$ is on the left side of $x=\tfrac12$. At $x=-\ell$, the two affine side heights are
	\[
	\widetilde h_\ell=\frac{i+p(\ell+1/2)}q\qquad\text{and}
	\qquad\widetilde g_\ell=\frac{j+p'(\ell+1/2)}{q'}.
	\]
	The slice on $x=-\ell$ cut by $\Delta_{i,j}$ is  $[\widetilde h_\ell,\widetilde g_\ell]$, and avoids integers exactly when \(\widetilde g_\ell<\lceil\widetilde h_\ell\rceil\),
	equivalently $j\le b_i^{\ell}$. This condition again holds automatically beyond the vertex.
	
	To combine the two directions, observe that the real thresholds satisfy
	\begin{align*}
		p'\Big(\ell+\frac12\Big)+q'\lfloor h_\ell\rfloor
		&\le c-\frac mq\Big(\ell+\frac12\Big)<c,\\
		-p'\Big(\ell+\frac12\Big)+q'\lceil\widetilde h_\ell\rceil
		&\ge c+\frac mq\Big(\ell+\frac12\Big)>c.
	\end{align*}
	It follows that \(a_i^{\ell}\le\lceil c\rceil\), and that 
	\(b_i^{\ell}\ge\lfloor c\rfloor\).
	Thus the upper bounds impose no additional restriction on integer labels when the third vertex is on the right side of $x=\tfrac12$, and the lower bounds impose none on integer labels when the third vertex is on the left side of $x=\tfrac12$. If $j=c$ is integral, all inequalities hold, consistently with the degenerate contribution. Hence the full criterion is
	\[
	\sup_{\ell\ge0}a_i^{\ell}\le j\le
	\inf_{\ell\ge0}b_i^{\ell}.
	\]
	
	The infinite extrema reduce to finite ones because
	\[
	a_i^{\ell+q}=a_i^{\ell}-m<a_i^{\ell}\qquad\text{and}
	\qquad b_i^{\ell+q}=b_i^{\ell}+m>b_i^{\ell}.
	\]
	 This proves the criterion.
	The symmetry identities follow directly from the formulas. Finally,
	\[
	 a_i\le\lceil c\rceil\le\lfloor c\rfloor+1
	\le b_i+1,
	\]
	which proves the non-negativity assertion.
\end{proof}

\subsection{The generating polynomial}\label{subsec:triangle-generating-polynomial}

We may introduce a generating polynomial which encodes the information about contributing Whitney triangles. 
As in the introduction, let  $\eta_{i,j}=\eta(\Delta_{i,j})\in\{0,1\}$ for the contribution indicator.

Let $ \mathbb F=\mathbb F_2$ denote the field with two elements. Compatible with the above identifications of index pairs, for commuting formal variables $\Xvar,\Yvar$, consider the algebra
\[
R=\mathbb F[\Xvar,\Yvar]/\langle \Xvar^{2q}\Yvar^{2q'}-1\rangle=\mathbb F[\Xvar,\Yvar]/\langle \Tvar^2-1\rangle
\qquad\text{where}\qquad \Tvar=\Xvar^q\Yvar^{q'}\in \mathbb F[\Xvar,\Yvar].
\]
Then $\Tvar$ may be regarded as an element of $R$ with $\Tvar^2=1$. Define
\begin{equation}\label{eq:general-triangle-polynomial}
	\mathcal P_{p/q,p'/q'}(\Xvar,\Yvar)=\sum_{(i,j)\in\Z^2/\langle(2q,2q')\rangle}
	\eta_{i,j}\Xvar^i\Yvar^j
	=\sum_{i=1}^{2q}\sum_{j\in\Zbb}\eta_{i,j}\Xvar^i\Yvar^j\in R.
\end{equation}
We usually abbreviate this polynomial to $\mathcal P$. Both variables are units in $R$: for example,
\[
\Xvar^{-1}=\Xvar^{2q-1}\Yvar^{2q'}\qquad\text{and}\qquad
\Yvar^{-1}=\Xvar^{2q}\Yvar^{2q'-1}.
\]
Thus negative exponents are well-defined. The sum is finite by the finite row intervals in Theorem~\ref{thm:general-row-interval}; for $m<0$ it is zero. Since equivalent labels give the same monomial, any choice of one representative of each contributing class gives the same element of $R$.

Translation by $(0,1)$ sends $(i,j)$ to $(i+q,j+q')$, preserves the curve families, and interchanges the two types of marked points. Since all marked points are avoided, it preserves the contribution indicator and \(\eta_{i+q,j+q'}=\eta_{i,j}\).

For $m>0$, this also follows from $a_{i+q}=a_i+q'$ and $b_{i+q}=b_i+q'$. Consequently, for all the slope pairs under consideration,
\begin{equation}\label{eq:triangle-polynomial-factor}
	\mathcal P=(1+\Tvar)\sum_{i=1}^{q}\sum_{j\in\Zbb}\eta_{i,j}\Xvar^i\Yvar^j
	=(1+\Tvar)\sum_{j=1}^{q'}\sum_{i\in\Zbb}\eta_{i,j}\Xvar^i\Yvar^j.
\end{equation}
The defining relation gives
\begin{equation}\label{eq:T-relations}
	\Tvar^2=1,\qquad (1+\Tvar)\Tvar^k=1+\Tvar\quad(k\in\Zbb)\qquad\text{and}
	\qquad (1+\Tvar)(1-\Tvar)=0.
\end{equation}
In particular, $1+\Tvar$ must not be canceled, since it is a zero-divisor.

We interpret fractions with denominators $\Xvar-1$ and $\Yvar-1$ in the localization
\[
 R[(\Xvar-1)^{-1},(\Yvar-1)^{-1}]
\]
of $R$.
The natural map from $R$ to this localization is injective. Indeed, the irreducible polynomials $\Xvar-1$ and $\Yvar-1$ do not divide $\Xvar^{2q}\Yvar^{2q'}-1$, as substitution of $\Xvar=1$ or $\Yvar=1$ shows. Each is relatively prime to the defining polynomial, hence is a non-zero-divisor in $R$. The formulas below identify elements of $R$, so these formal quotients are uniquely determined there; this does not assert that the denominators are units in $R$.

By \eqref{eq:triangle-polynomial-factor} and Theorem~\ref{thm:general-row-interval},
the general finite formula is
\begin{equation}\label{eq:general-generating-floor-formula}
		\mathcal P=(1+\Tvar)\sum_{i=1}^{q} \Xvar^i\Big(\Yvar^{a_i}+\Yvar^{a_i+1}+\cdots+\Yvar^{b_i}\Big)
		=\frac{1+\Tvar}{(\Yvar-1)}\Big(\sum_{i=1}^{q}
		\Big(\Xvar^i\Yvar^{b_i+1}-\Xvar^i\Yvar^{a_i}\Big)\Big).
\end{equation}
Here a row with $b_i=a_i-1$ contributes zero.

\subsection{An integer between the slopes}\label{subsec:integer-between-slopes}

The extrema $a_i$ and $b_i$ simplify when the closed interval between the slopes contains an integer, i.e. when
\begin{equation}\label{eq:integer-between-slopes}
	\Big\lceil\frac{p'}{q'}\Big\rceil
	\le\Big\lfloor\frac pq\Big\rfloor.
\end{equation}
Write the first-column endpoints explicitly as
\[
a_i^{0}=\frac{p'+1}{2}+q'\Big\lfloor\frac{2i-p}{2q}\Big\rfloor\qquad\text{and}
\qquad b_i^{0}=q'\Big\lceil\frac{2i+p}{2q}\Big\rceil-\frac{p'+1}{2},
\]
and put
\[
\Mvar=\Xvar^{(p+1)/2}\Yvar^{(p'+1)/2}\qquad\text{and}\qquad
\Nvar=\Xvar^{(1-p)/2}\Yvar^{(1-p')/2}=\frac{\Xvar\Yvar}{\Mvar}.
\]

\begin{thm}\label{thm:integer-between-slopes}\label{cor:thm-when-m-is-large}\label{cor:separated-generating-polynomial}
	Assume $m>0$. The following are equivalent:
	\begin{enumerate}
		\item The integer-between-slopes condition \eqref{eq:integer-between-slopes} holds.
		\item For every $i\in\Zbb$, $a_i=a_i^{0}$ and $b_i=b_i^{0}$.
		\item The generating polynomial is given by
		\begin{equation}\label{eq:separated-polynomial-symmetric}
			\mathcal P
				=\frac{1+\Tvar}{(\Xvar-1)(\Yvar-1)}\bigl(\Mvar(1-\Xvar^q)+\Nvar(1-\Yvar^{q'})\bigr).
		\end{equation}
	\end{enumerate}
	Under these conditions,  there are $2m+2qq'$ contributing classes. 
\end{thm}
\begin{proof}
	We first establish the endpoint criterion, then compute $\mathcal P_0$, where $\mathcal P_\ell$ is defined by 
	\begin{align*}
		\mathcal P_\ell=(1+\Tvar)\sum_{i=1}^{q} \Xvar^i\Big(\Yvar^{a_i^{\ell}}+\cdots+\Yvar^{b_i^{\ell}}\Big)
		=\frac{1+\Tvar}{(\Yvar-1)}\Big(\sum_{i=1}^{q}
		\Big(\Xvar^i\Yvar^{b_i^{\ell}+1}-\Xvar^i\Yvar^{a_i^{\ell}}\Big)\Big).
	\end{align*}

	Suppose \eqref{eq:integer-between-slopes} holds and put $a=\lfloor\tfrac pq\rfloor$. Then $\tfrac{p'}{q'}\le a\le \tfrac pq$. For every $\ell\ge0$,
	\begin{align*}
		a_i^{\ell}-a_i^{0}
		&=p'\ell+q'\Big(
		\Big\lfloor\frac{2i-p}{2q}-\frac{p\ell}{q}\Big\rfloor
		-\Big\lfloor\frac{2i-p}{2q}\Big\rfloor\Big)
		\le\ell(p'-aq')\le0.
	\end{align*}
	Therefore $ a_i=a_i^{0}$, and $ b_i=b_i^{0}$ follows by the reflection identity. This proves (1)$\Rightarrow$(2).
	
	Conversely, suppose the interval between the slopes contains no integer. Write
	\[
	p=aq+r,\qquad a=\Big\lfloor \frac pq\Big\rfloor,\qquad 1\le r\le q-1.
	\]
	Here $q\ge2$ and $p'>aq'$. Let $i\in\{1,\ldots,q\}$ denote the residue of $\tfrac{p-1}2$ modulo $q$.
	The fractional part of $\tfrac{2i-p}{2q}$ is $1-\tfrac1{2q}$, which is larger than $\tfrac rq$. Hence
	\begin{align*}
	\Big\lfloor\frac{2i-p}{2q}-\frac pq\Big\rfloor
	=\Big\lfloor\frac{2i-p}{2q}\Big\rfloor-a,
	\qquad\text{and}\qquad
	a_i^{1}-a_i^{0}=p'-aq'>0.
	\end{align*}
	Since $1<q$, this column occurs in the finite maximum defining $ a_i$. Thus (2) fails. We have proved (1)$\Leftrightarrow$(2).
	
	We now compute the first-column polynomial $\mathcal P_0$. As before, we assume that $m>0$. Note that
	
	\begin{equation}\label{eq:polynomial-endpoint-difference}
		(\Xvar-1)(\Yvar-1)\mathcal P_0
		=(1+\Tvar)(\Xvar-1)\sum_{i=1}^{q}\Xvar^i(\Yvar^{b_i^{0}+1}-\Yvar^{a_i^{0}}).
	\end{equation}
	Choose the unique integers $s,t$ with \(p+2t+1=2qs\) and \(0\le t\le q-1\).
	For $1\le i\le q$,
	\[
	\Big\lceil\frac{2i+p}{2q}\Big\rceil
	=\begin{cases}s,&i\le t,\\s+1,&i>t,\end{cases}
	\qquad
	\Big\lfloor\frac{2i-p}{2q}\Big\rfloor
	=\begin{cases}-s,&i<q-t,\\1-s,&i\ge q-t.\end{cases}
	\]
If we set \(a=sq'+\tfrac{1-p'}2\) and \(b=-sq'+\frac{p'+1}2\),
	the two sums are
	\begin{align*}
		(\Xvar-1)\sum_{i=1}^{q}\Xvar^i\Yvar^{b_i^{0}+1}
		&=(\Xvar-1)\Yvar^a\Big(\sum_{i=1}^{t}\Xvar^i
		+\Yvar^{q'}\sum_{i=t+1}^{q}\Xvar^i\Big)\\
		&=\Yvar^a\bigl(\Xvar^{t+1}-\Xvar+\Yvar^{q'}\Xvar^{q+1}-\Yvar^{q'}\Xvar^{t+1}\bigr)
		\qquad\text{and}\\
		(\Xvar-1)\sum_{i=1}^{q}\Xvar^i\Yvar^{a_i^{0}}
		&=(\Xvar-1)\Yvar^b\Big(\sum_{i=1}^{q-t-1}\Xvar^i
		+\Yvar^{q'}\sum_{i=q-t}^{q}\Xvar^i\Big)\\
		&=\Yvar^b\bigl(\Xvar^{q-t}-\Xvar+\Yvar^{q'}\Xvar^{q+1}-\Yvar^{q'}\Xvar^{q-t}\bigr).
	\end{align*}
	Subtracting and using $\Tvar=\Xvar^q\Yvar^{q'}$, we find
	\[\mathcal P_0=\frac{(\Tvar+1)}{(\Xvar-1)(\Yvar-1)}\Big(
	\Xvar(\Tvar-1)(\Yvar^a-\Yvar^b)
	+(1-\Yvar^{q'})(\Xvar^{t+1}\Yvar^a-\Xvar^{q-t}\Yvar^b)\Big).
	\]
	The first term inside the parentheses vanishes after multiplication by $1+\Tvar$.
	Since $qs-t=\tfrac{p+1}2$,
	\begin{align*}
		&(1+\Tvar)\Xvar^{t+1}\Yvar^a
		=(1+\Tvar)\Xvar^{t+1-qs}\Yvar^{\frac{1-p'}2}=(1+\Tvar)\Nvar,\\
		&(1+\Tvar)\Xvar^{q-t}\Yvar^b
		=(1+\Tvar)\Xvar^{q-t+qs}\Yvar^{\frac{p'+1}2}=(1+\Tvar)\Xvar^q\Mvar\\
		\Rightarrow\qquad&\mathcal P_0
		=\frac{(\Tvar+1)}{(\Xvar-1)(\Yvar-1)}(\Yvar^{q'}-1)(\Xvar^q\Mvar-\Nvar)\\
		&\quad=\frac{(\Tvar+1)}{(\Xvar-1)(\Yvar-1)}
		(\Mvar+\Nvar-\Xvar^q\Mvar-\Yvar^{q'}\Nvar),
	\end{align*}
	where the last equality follows since $(1+\Tvar)\Xvar^q\Yvar^{q'}\Mvar=(1+\Tvar)\Mvar$.
	 It follows that (3) is equivalent to $\mathcal P=\mathcal P_0$.
	Now, if  (2) holds,  
	$\mathcal P=\mathcal P_0$. If (2) fails, some lower or upper tail is nonempty. The monomials $\Xvar^i\Yvar^j$ with $1\le i\le2q$ form a basis of the group algebra $\mathbb F[\mathbb Z^2/\langle(2q,2q')\rangle]=R$. Distinct removed labels therefore cannot cancel, even over $\mathbb F_2$, and $\mathcal P\ne\mathcal P_0$. This proves (2)$\Leftrightarrow$(3), as well as uniqueness in (3).
	
	Under these equivalent conditions, the contribution criterion follows from Theorem~\ref{thm:general-row-interval}, and the endpoint symmetries follow directly from their formulas. For the class count,
	\[
	b_i^{0}-a_i^{0}+1=q'\Big(
	\Big\lceil\frac{2i+p}{2q}\Big\rceil
	-\Big\lfloor\frac{2i-p}{2q}\Big\rfloor\Big)-p'.
	\]
	Using
	\[
	\sum_{i=1}^{q}\Big\lfloor\frac{i+c}{q}\Big\rfloor=c+1
	\qquad(c\in\Zbb)
	\]
	with $c=(p-1)/2$ and $c=-(p+1)/2$, the parenthesized difference sums to $p+q$ over $1\le i\le q$. It is $q$-periodic, so the total count is
	\(
	q'(2p+2q)-2qp'=2m+2qq'
	\).
\end{proof}

The formula (\ref{eq:separated-polynomial-symmetric}) also displays 
the symmetry
\begin{equation}\label{eq:polynomial-slope-symmetry}
	\mathcal P_{p/q,p'/q'}(\Xvar,\Yvar)
	=\mathcal P_{-p'/q',-p/q}(\Yvar,\Xvar).
\end{equation}
This substitution preserves $m$, $\Tvar$, and the integer-between-slopes condition, while interchanging $\Mvar,\Nvar$ and $\Xvar^q,\Yvar^{q'}$. More generally, the aforementioned symmetry follows geometrically by reflecting $(x,y)\mapsto(1-x,y)$ and interchanging $\beta,\gamma$, which sends $(i,j)$ to $(j,i)$ and preserves the contribution criterion. Thus \eqref{eq:polynomial-slope-symmetry} also holds without the integer-between-slopes hypothesis.

\section{Tool Use}
GPT-6 assisted the authors with several aspects of this paper. The bi-grading theorem (Theorem~\ref{thm:line-maslov-bigrading}) was formulated and proved through correspondence with GPT-6. The authors initially conjectured that Theorem~\ref{thm:2} might hold for all pairs of rational slopes. GPT-6 provided a counterexample to this conjecture and also helped establish the equivalence of the three conditions in Theorem~\ref{thm:integer-between-slopes}. The latter result reformulates an earlier version of the theorem to make the equivalence of these conditions explicit.

GPT-6 also assisted with proofreading, grammatical corrections, and improvements to the exposition and overall flow of the paper. The authors have reviewed and verified the resulting mathematical arguments and take full responsibility for the contents of the paper.

\bibliographystyle{hamsalpha}
\newcommand{\etalchar}[1]{$^{#1}$}

\end{document}